\documentclass{amsart}

\usepackage[all]{xy}

\usepackage{latexsym}
\usepackage{amsmath}
\usepackage{amssymb}

\def\eu{\mathfrak}
\def\ma{\mathbb}

\def\fin{\hfill\qed\bigskip}
\def\p{{\eu p}_{\infty}}

\newcommand{\Aut}{\operatorname{Aut}}
\newcommand{\PSL}{\operatorname{PSL}}
\newcommand{\id}{\operatorname{Id}}
\newcommand{\Gal}{\operatorname{Gal}}
\newcommand{\con}{\operatorname{con}}

\newcommand{\xbinom}{\genfrac[]{0.5pt}0}

\numberwithin{equation}{section}
\newtheorem{theorem}{Theorem}[section]
\newtheorem{remark}[theorem]{Remark}

\title[Congruence Function Fields with Class Number One]
 {Congruence Function Fields with Class Number One}
 
\dedicatory{To the memory of Prof. Manohar L. Madan}

\author[M. Rzedowski]{Martha Rzedowski--Calder\'on}
\address{Departamento de Control Autom\'atico\\
Centro de Investigaci\'on y de Estudios Avanzados del I.P.N.}
\email{mrzedowski@ctrl.cinvestav.mx}

\author[G. Villa]
{Gabriel Villa--Salvador}
\address{
Departamento de Control Autom\'atico\\
Centro de Investigaci\'on y de Estudios Avanzados del I.P.N.}
\email{gvillasalvador@gmail.com, gvilla@ctrl.cinvestav.mx}



\date{December 15, 2014}

\begin{document}

\begin{abstract}
We prove that there exists, up to isomorphism,
exactly one function field over the
finite field of two elements of class number one and genus four.
This result, together with the ones of MacRae, Madan, Leitzel,
Queen and Stirpe, establishes that there exist eight non-isomorphic
congruence function fields of genus larger than zero and class
number one.
\end{abstract}

\maketitle

\section{Introduction}\label{S1}

Let $K$ be a congruence function field with exact field of constants
${\ma F}_q$, the finite field of $q$ elements. Consider the class
group of divisors of degree zero of $K$: $C_{0K}$. It is a finite abelian
group with $h_K$ elements, $h_K$ is called {\em the class number of
$K$}. When $K$ is a function field of genus $0$, we have $h_K=1$. Thus,
we consider $K$ of genus $g_K\geq 1$. When $q\geq 5$ and
$g_K\geq 1$ we have $h_K>1$. In \cite{Mac71} R. MacRae found
all the congruence function fields with class number one in the
particular case that $K$ is a quadratic extension of the rational function
field $k={\ma F}_q(T)$ and $K$ contains a prime divisor of degree
one. He proved that there are four quadratic fields with class 
number one which have a prime of degree one.
M. Madan and C. Queen continued
the study of this problem in \cite{MaQu72}. They
showed that if $q=2$ and $g_K>4$, or $q=3$ and $g_K>2$
then $h_K\neq 1$. Finally, they proved that except for the case $q=2$
and $g_K=4$ there exist exactly seven congruence function fields
with class number one and genus larger than zero. The case $q=2$,
$g_K=4$ was not settled.

In \cite{LeMaQu75} J. Leitzel, M. Madan and C. Queen considered the case
$q=2$ and $g_K=4$ and claimed
that there is no field of class number one over the finite field of two
elements and genus four. However C. Stirpe \cite{Sti2013} found a 
counterexample to this claim. The example runs as follows.
Let ${\eu m}$ be the place associated to the irreducible
polynomial $T^4+T+1\in {\ma F}_2[T]$ and let 
$S$ be the place associated to the irreducible polynomial $T^7+T^4+1$.
Let $K^{\eu m}_S$ be the ray class field of conductor ${\eu m}$
and such that $S$ splits in $K^{\eu m}_S/{\ma F}_2(T)$. Stirpe established
that the subfield of degree five over ${\ma F}_2(T)$ satisfies that
$h_K=1$ and $g_K=4$. Furthermore,
Stirpe claims that $T^7+T^4+1$ is not unique. For instance, he
remarks that we may take $S_1$ to be the place associated to $T^7+
T^3+1$ and the unique subfield $K_1$ of $K^{\eu m}_{S_1}$ of 
degree five over ${\ma F}_2(T)$ also satisfies that $h_{K_1}=1$
and $g_{K_1}=4$.

In \cite{StiMer2014}, P. Mercuri and C. Stirpe proved that the 
two fields found by Stirpe in \cite{Sti2013} are in
fact isomorphic. Furthermore, they show that, up to isomorphism,
there is only one field of genus $4$ and class number one.
This result together with the results of Madan, Leitzel, Queen and
Stirpe, shows that, up to isomorphism,
there are exactly eight congruence function fields
$K$ of genus larger than zero and class number one.

In this paper we present another proof that, up to isomorphism,
there is only one field of genus four and class number one. 
We do not use the examples found by Stirpe in \cite{Sti2013}.
Our approach uses the theory of cyclotomic function 
fields of Carlitz--Hayes. First, we
consider a field $K$ over ${\ma F}_2$ such that $h_K=1$ and 
$g_K=4$. We show that $K$ has a unique rational function subfield
$k:={\ma F}_2(T)$ such that $[K:k]=5$ and that the extension
$K/k$ is cyclic. There is only one prime of $k$ ramified in
$K$ and this place is of degree four. From the result of Madan and
Queen \cite{MaQu72} that states that a function field $K$ over
${\ma F}_2$ satisfies that $h_K=1$ and $g_K=4$ if and only if
$N_1=N_2=N_3=0$ and $N_4=1$ where $N_i$ denotes the number
of prime divisors in $K$ of degree $i$, we deduce that, up to
isomorphism, necessarily $K\subseteq k(\Lambda_M){\ma F}_{2^5}$
where $M=T^4+T+1$ and $k(\Lambda_M)$ is the cyclotomic function
field corresponding to the Carlitz module $\Lambda_M$.
Finally, we prove that there are precisely two fields $K$ over $k$  contained
in $k(\Lambda_M){\ma F}_{2^5}$ such that $N_1=N_2=N_3=0$ and
$N_4=1$. In one of them $T^7+T^4+1$ splits and $T^7+T^3+1$
is inert and in the other $T^7+T^4+1$ is inert and $T^7+T^3+1$
splits. Both fields are isomorphic.

One of the key facts in our proof is that if $[K:k]
=5$ and ${\eu p}$ is the divisor of degree four in $K$, then the different
of $K/k$ is ${\eu p}^4$ and ${\eu p}$ is totally ramified. This
was proved by Mercuri and Stirpe in \cite{StiMer2014}.

\section{The field $K$}\label{S2}

Let $K$ be a congruence function field with
exact field of constants the finite field of $q$ elements ${\ma F}_q$.
Let $N_i$ denote the number of prime divisors of degree $i$ in $K$,
$i\geq 1$. Let $A_i$ be the number of integral divisors in $K$ of
degree $i$, $i\geq 0$. The genus of $K$ will be denoted by $g$
and the class number of $K$ by $h$. Let $k={\ma F}_q(T)$ be a
rational congruence function field and let $R_T={\ma F}_q[T]$ be
its ring of integers. $\p$ will denote the pole divisor of $T$
in $k$. For the standard results on congruence function fields
and cyclotomic function fields we refer to \cite{Vil2006}.

For any divisor ${\eu q}$ in $K$ we denote by $d_K({\eu q})$
its degree. If $P_K(u)=a_0+a_1u+\cdots+a_{2g}u^{2g}$ is the
numerator of the zeta function of $K$, where $u=q^{-s}$, we have
the following relations (\cite[Theorems 6.3.5 and 6.4.1]{Vil2006})
\begin{gather}
a_0=1,\quad a_{2g}=q^g, \quad a_{2g-i}=a_iq^{g-i}, \quad 0\leq 
i\leq 2g,\nonumber\\
a_i=A_i-(q+1)A_{i-1}+qA_{i-2}, \quad 0\leq i \leq 2g,
\quad \text{with}\quad A_{-1}=
A_{-2}=0,\label{Eq-1}\\
P_K(1)=h, \quad A_n=h\Big(\frac{q^{n-g+1}-1}{q-1}\Big)\quad
\text{for}\quad n>2g-2.\nonumber
\end{gather}

From now on, $K$ will denote a field over ${\ma F}_2$ such that $g=4$
and $h=1$. This condition is equivalent
to $N_1=N_2=N_3=0$ and $N_4=1$ (\cite[Theorem 2 (v)]{MaQu72}).
From that paper we know that the numerator of the zeta function
of $K$ is $P_K(u)=1-3u+2u^2+u^4+8u^6-24u^7+16u^8$. Let ${\eu p}$ denote
the only prime divisor of degree four in $K$.

In this case, from (\ref{Eq-1}) we obtain that $A_0=1$, 
$A_i=N_i=0$, $0\leq i\leq 3$, $
A_4=N_4=1$ and $A_5=N_5=3$. Let ${\eu C}_i$, $1\leq i\leq 3$ be the
three places of degree five in $K$. Therefore $\ell({\eu C}_1^{-1})=
2$, $L({\eu C}_1^{-1})=\{0,1,T,T+1\}$ where $(T)_K=\frac{{\eu C}_2}
{{\eu C}_1}$ and $(T+1)_K=\frac{{\eu C}_3}{{\eu C}_1}$ where
$(y)_K$ denotes the divisor in $K$ of $y\in K^{\ast}$. We have
$[K:k]=5$. Since $L({\eu p}^{-1})={\ma F}_2$,
it follows that the minimal $n$ such that there exists $y\in K$
with $[K:{\ma F}_2(y)]=n$ is $n=5$ and that $k$ is unique
satisfying this property.

\begin{remark}\label{R0.-1}{\rm{
Every proper subfield ${\ma F}_2\subsetneqq E\subsetneqq
K$ such that $K/E$ is separable,
is of genus $0$. Indeed, for any finite subextension $E$ of
$K$, the differential exponent of every prime appearing in the
different ${\eu D}_{K/E}$ of the extension is greater than or equal to
$2$ and since the minimum degree of a prime in $K$ is $4$, the
degree of $d$ of ${\eu D}_{K/E}$ is greater than of equal to $8$ except
in the case that $K/E$ is unramified. From the 
Riemann-Hurwitz formula, if $g_E\geq 1$ and $K/E$ ramified,
we obtain
\[
6=2g_K-2=[K:E](2g_E-2)+d\geq d\geq 8.
\]
Thus, if $g_E\geq 1$, $K/E$ is unramified, $[K:E]=3$, and $g_E=2$. 
If $K/E$ is normal, let ${\eu t}={\eu p}\cap E$. Then since $[K:E]=3$
is relatively prime to $\deg_K {\eu p}=4$, it follows that $\eu t$ 
decomposes fully in $K/E$ and in particular $K$ would contain at
least $3$ primes of degree four. Therefore $K/E$ is non-normal.
Let $\tilde{K}$ be the Galois closure of $K/E$. Then $[\tilde{K}:K]
=2$ and since $\tilde{K}/K$ is unramified, it follows that $\tilde{K}
=K{\ma F}_4$. We have that $\tilde{K}/E{\ma F}_4$ is a normal
extension of degree $3$. Since $\deg_K {\eu p}=4$ and $\tilde{K}
/K$ is an extension of constants of degree $2$, we obtain that
$\eu p$ decomposes into two primes of degree $2$ in $\tilde{K}$
(see \cite[Theorem 6.2.1]{Vil2006}). Thus $\tilde{K}$ has exactly
two primes of degree $2$. Let $\tilde{\eu p}$ one of them and let
$\tilde{\eu t}=\tilde{\eu p}\cap E{\ma F}_4$. As above we obtain that
$\tilde{\eu t}$ decomposes fully in $\tilde{K}/E{\ma F}_4$ and
in particular we have at least three primes in $\tilde{K}$ of 
degree $2$. This contradiction shows that $g_E=0$.
}}
\end{remark}

\begin{remark}\label{R0.1}{\rm{
Let $\theta\in{\mathcal G}:=\Aut_{{\ma F}_2} K$. 
Since $\theta$ permutes
the three divisors ${\eu C}_i$, $1\leq i\leq 3$, we have that
$\theta|_{k}\in \Aut_{{\ma F}_2} k\cong \PSL(2,{\ma F}_2)
\cong S_3$ where $k={\ma F}_2(T)$ and $S_3$ is the symmetric
group in three elements. Therefore $K^{{\mathcal G}}\supseteq
k^{S_3}$. Therefore $|{\mathcal G}|$ divides $30$.
If $5$ divides $|{\mathcal G}|$, then
the field fixed by an element of order $5$ of 
${\mathcal G}$  is necessarily $k$ and $K/k$ is normal.
If $K/k$ is not normal then ${\mathcal G}$ is trivial since otherwise
for each non--trivial subgroup of ${\mathcal G}$, the fixed field is
of genus $0$ but one of them
is of degree less than five. This contradicts that five is the minimum
degree of a proper subfield of $K$. Therefore, we have that
$K/k$ is normal if and only if $\Aut_{{\ma F}_2} K\neq\{\id\}$.
}}
\end{remark}

One of the key facts to prove the uniqueness of $K$
is the following theorem.

\begin{theorem}\label{T2.1} The extension $K/k$ is
normal.
\end{theorem}

\proof Stirpe and Mercuri \cite{StiMer2014} proved that ${\eu p}$
is fully ramified in $K/k$ and in particular
${\eu D}_{K/k}={\eu p}^4$, where ${\eu D}_{K/k}$
denotes the different of the extension $K/k$. 

Assume that $K/k$ is not normal. Let
$\tilde{K}$ be the Galois closure of $K/k$, $G:=\Gal(\tilde{K}/k)$ and
$H:=\Gal(\tilde{K}/K)$. Then $G$ is a transitive subgroup of $S_5$,
the symmetric group in five elements and $H$ is a subgroup of $S_4$.
The field of constants of $\tilde{K}$ is ${\ma F}_2$ because otherwise,
since the primes of degree one are inert in $K/k$, we would have
an element in $G$ of order $5 r$ with $r\geq 2$ contrary to the fact
that the elements in $S_5$ are of order less than or equal to six.

From Abhyankar Lemma we obtain that
$\tilde{K}/K$ is unramified. Let $H_1$ be a proper
normal subgroup of $H$ such that $H/H_1$ is abelian. Then
we obtain a non--trivial unramified 
abelian extension of $K$ and since the class number of $K$ is one,
this extension would be a constant extension. This contradiction
proves that $K/k$ is normal.  \fin

We have 
${\eu D}_{K/k}={\eu p}^4$. Since $N_1=N_2=N_3=0$
and $N_4=1$, we obtain that all prime divisors of $k$
of degree less than or equal to four, except
for one of degree four, are
inert in $K/k$ and one prime divisor of degree four is
ramified.

In $k$ we have three prime divisors of degree four, namely
the ones corresponding to $T^4+T+1$, $T^4+T^3+1$ and
$T^4+T^3+T^2+ T+1$.

\begin{remark}\label{R2.2} {\rm{
We may assume without loss of generality
that the ramified prime of degree four ${\eu m}$ 
is the place corresponding
to $M=T^4+T+1$ for if ${\eu m}_1$ is the place 
corresponding to $T^4+T^3+1$
(resp. $T^4+T^3+T^2+T+1$), then $\sigma\colon k\to
k$ given by $\sigma(T)=\frac{1}{T}$ (resp. $\sigma(T)=\frac{1}
{T+1}$) satisfies $\sigma(T^4+T+1)=\frac{T^4+T^3+1}{T^4}$ (resp.
$\sigma(T^4+T+1)=\frac{T^4+T^3+T^2+T+1}{(T+1)^4}$) 
so that $\sigma({\eu m})=
{\eu m}_1$ and extending $\sigma$ to $\tilde{\sigma}\colon K\to
\overline{k}$ we obtain $\tilde{\sigma}(K)\cong K$ and
$\sigma(k)=k$. 
In $\tilde{\sigma}(K)/k$, the prime
${\eu m}_1$ is the ramified one.}}
\end{remark}

The extension $K/k$ is a cyclic extension such
that all the primes of degree one, two and 
three in $k$ ($\{\p,T,T+1,T^2+T+1,T^3+T+1,T^3+T^2+1\}$) and the 
primes of degree four associated to $T^4+T^3+1$ and $T^4+T^3+T^2
+T+1$ are inert. The prime ${\eu m}$ associated to
$T^4+T+1$ is ramified.

Since $\p$ is unramified in $K/k$, in fact $\con_{k/K}\p={\eu C}_1$, and
${\eu m}$ is the only ramified prime in $K$ and it is tamely ramified
we have that $K\subseteq k(\Lambda_M){\ma F}_{2^5}$ (see
\cite[Proposition 3.4]{MaRzVi2013}) and $[K:k]=5$.

The key step for the main result of this paper is the following theorem.

\begin{theorem}\label{T2.3}
Up to isomorphism, there exists only one field $K$ with $k\subseteq K
\subseteq k(\Lambda_M){\ma F}_{2^5}$ such that $g=4$ and $h=1$.
\end{theorem}

\proof
To start, let ${\eu t}$ be any prime divisor of $k$ such that ${\eu t}
\neq \p, {\eu m}$, and let $P\in R_T:={\ma F}_2[T]$ be the monic
irreducible polynomial associated to ${\eu t}$. Then the Frobenius
map $\varphi_P$ of $P$ in the extension $k(\Lambda_M)/k$ is given
by $\varphi_P(\lambda)=\lambda^P$ where $\lambda$ is a generator
of $\Lambda_M$ (see \cite[Theorem 12.5.1]{Vil2006}). In particular
for ${\eu t}\neq \p, {\eu m}$ we have that the decomposition group
of ${\eu t}$ is $D_P=\langle \varphi_P\rangle$ and $|D_P|=
o(P\bmod M)$.

We have that $G_M:=\Gal(k(\Lambda_M)/k)\cong C_{15}$, the cyclic
group of $15$ elements, and let $L$ be the subfield of $k(\Lambda_M)$
such that $[L:k]=5$.
\[
\xymatrix{k(\Lambda_M)\ar@{-}[r]\ar@{-}[d]_3^{\langle\tau^5\rangle}&
k(\Lambda_M){\ma F}_{2^5}\ar@{-}[d]\\
L\ar@{-}[r]\ar@{-}[d]_5^{\langle\tau^3\rangle}&L{\ma F}_{2^5}\ar@{-}[d]\\
k\ar@{-}[r]&k{\ma F}_{2^5}}
\]

From the isomorphism $G_M\cong \big(R_T/(M)\big)^{\ast}$ we have
that $\tau$, given by $\tau(\lambda)=\lambda^T$, is a generator
of $G_M$. Therefore $\Gal(L/k)\cong\langle \tau\bmod 
\langle\tau^5\rangle\rangle
\cong \langle \tau^3\rangle$.

Note that $P$ is inert in $L/k$ if and only if $o(\varphi_P)\in\{5,15\}$.
Direct computations give

\begin{gather}
T^1\equiv T\bmod M,\quad T^2\equiv T^2\bmod M,\quad T^3\equiv
T^3\bmod M,\quad T^4\equiv T+1\bmod M,\nonumber\\
T^5\equiv T^2+T\bmod M,\quad T^6\equiv T^3+T^2\bmod M,
\quad T^7\equiv T^3+T+1\bmod M,\nonumber\\
T^8\equiv T^2+1\bmod M,\quad T^9\equiv T^3+T\bmod M,\quad
T^{10}\equiv T^2+T+1\bmod M,\nonumber\\
T^{11}\equiv T^3+T^2+T\bmod M,\quad T^{12}\equiv T^3+T^2+
T+1\bmod M,\label{Eq2}\\
T^{13}\equiv T^3+T^2+1\bmod M,\quad \text{and}\quad
T^{14}\equiv T^3+1\bmod M,\nonumber\\
\text{and}\quad T^4+T^3+1\equiv T^3+T\bmod M, \quad 
T^4+T^3+T^2+T+1\equiv T^3+T^2\bmod M.\nonumber
\end{gather}

From (\ref{Eq2}) we may compute the order of $\varphi_P$:
\begin{gather*}
o(\varphi_T)=15,\quad o(\varphi_{T+1})=15,\quad o(\varphi_{T^2+T+1})
=3,\quad o(\varphi_{T^3+T^2+1})=15,\\
o(\varphi_{T^3+T+1})=15,\quad o(\varphi_{T^4+T^3+1})=5, \quad o(\varphi_{
T^4+T^3+T^2+T+1})=5,
\end{gather*}
and we also have that $\p$ is fully decomposed in $k(\Lambda_M)/k$
(see \cite[Theorem 12.4.6]{Vil2006}), that is $o(\varphi_{\p})=1$
where $\varphi_{\p}$ denotes the Frobenius of $\p$ in $k(\Lambda_M)
/k$. Therefore, the decomposition groups of $P$ in 
$k(\Lambda_M)/k$ are given by
\begin{gather}
D_T=D_{T+1}=D_{T^3+T^2+1}=D_{T^3+T+1}=G_M=\langle\tau\rangle,\nonumber\\
D_{\p}=\{\id\},\quad D_{T^2+T+1}=\langle\tau^5\rangle,\label{Eq0}\\
D_{T^4+T^3+1}=D_{T^4+T^3+T^2+T+1}=\langle\tau^3\rangle.\nonumber
\end{gather}
In particular $\p$ and $T^2+T+1$ are decomposed in $L/k$ and $T,T+1,
T^3+T^2+1,T^3+T+1, T^4+T^3+1$ and $T^4+T^3+T^2+T+1$ are
inert in $L/k$.

Now, in the extension of constants $k_5:=k{\ma F}_{2^5}$
over $k$, all the primes
of degree $i$, $1\leq i\leq 4$, are inert (\cite[Theorem 6.2.1]{Vil2006}).
We have that $\Gal(k_5/k)=\langle\chi\rangle$ where $\chi$ is induced
by the Frobenius map of the extension ${\ma F}_{2^5}/{\ma F}_2$. More
precisely, if $Q(T)=\sum_{i=0}^d a_iT^i\in {\ma F}_{2^5}[T]$, then $\chi(Q(T))=
\sum_{i=0}^d a_i^2T^i$.

Let $P(T)\in R_T$ be a prime of degree $i$, $0\leq i\leq 4$. Then the
residue fields in $k_5/k$ are isomorphic to ${\ma F}_{2^{5i}}/
{\ma F}_{2^i}$ and the Frobenius map $\delta$ of ${\ma F}_{2^{5i}}/
{\ma F}_{2^i}$ is given by $\delta(\alpha)=\alpha^{2^i}$ for $\alpha
\in {\ma F}_{2^{5i}}$. Therefore the Frobenius map of $P(T)$
in $k_5/k$ corresponds
to $\langle\chi^i\rangle \in \Gal(k_5/k)$.

To find the Frobenius map of an arbitrary $P\in R_T$ in the extensions
$Lk_5$ and $k(\Lambda_M)k_5$ we consider the following general
situation.
Let $E/F, J/F$ be Galois extensions of global or local fields such that
$E\cap J=F$. Let $S:=EJ$. 
\[
\xymatrix{E\ar@{-}[d]\ar@{-}[r] &EJ=S\ar@{-}[d]\\
F\ar@{-}[r]&J}
\]

We have the isomorphism
\begin{align*}
\Phi\colon\Gal(S/F)&\to \Gal(E/F)\times \Gal(J/F)\\
\Phi(\theta)&=\big(\theta|_E,\theta|_J),
\end{align*}
and the inverse of $\Phi$ is given by
\begin{align*}
\Psi\colon \Gal(E/F)\times\Gal(J/F)&\to\Gal(S/F)\\
\Psi(\alpha,\beta)&=\tilde{\alpha}\tilde{\beta},
\end{align*}
where $\tilde{\alpha}\colon S\to S$ and $\tilde{\beta}\colon S\to S$ are
defined, for $z=\sum_{i=1}^t x_iy_i\in S$ with $x_i\in E$ and $y_i\in J$,
by
\begin{align*}
\tilde{\alpha}\big(\sum_{i=1}^t x_iy_i\big)&=\sum_{i=1}^t\alpha(x_i)y_i
\intertext{and}
\tilde{\beta}\big(\sum_{i=1}^t x_iy_i\big)&=\sum_{i=1}^t x_i\beta(y_i).
\end{align*}

Let ${\mathcal P}$ be a prime in $F$, ${\eu P}$ be a prime in $S$ above
${\mathcal P}$ and let ${\eu q}:={\eu P}\cap J$ and ${\eu t}:={\eu P}\cap E$.
Assume that ${\mathcal P}$ is unramified in $S/F$. Let $\xbinom{S/F}
{{\eu P}}\in \Gal(S/F)$ be the Frobenius map of ${\eu P}/{\mathcal P}$. Then
$\xbinom{S/F}{{\eu P}}\Big|_E=\xbinom{E/F}{{\eu t}}$ and
$\xbinom{S/F}{{\eu P}}\Big|_J=\xbinom{J/F}{{\eu q}}$. Therefore
\begin{gather}\label{Eq1}
\xbinom{S/F}{{\eu P}}=\widetilde{\xbinom{E/F}{{\eu t}}}
\widetilde{\xbinom{J/F}{{\eu q}}}.
\end{gather}

We will apply formula (\ref{Eq1}) to our case $(F=)k={\ma F}_2(T)$,
$(E=)L$, $(J=)k_5$ and $(S=)Lk_5$. There exist exactly
four extensions $R_j$, $1\leq j\leq 4$
of degree five over $k$ contained in $L_5:=Lk_5$ other than $L$
and $k_5$. 
The fields $K$ we are looking are,
if any, among the fields $R_j$ such that all the primes $P(T)$
in $k$ of degree less than or
equal to four other than $M$ are inert in $K/k$.
Note that since
$k_5/k$ is unramified and the only ramified prime in $L/k$ is 
${\eu m}$, the only ramified prime in each $R_j$ is ${\eu m}$.

We have $\Gal(L_5/k)\cong C_5\times C_5$ and the decomposition group
of any unramified prime is cyclic since the characteristic is $2\neq 5$.
Thus, any prime of degree $i$ with $i\leq 4$ other than ${\eu m}$
is decomposed in exactly one field among $L$, $R_j$, $1\leq j\leq 4$, namely,
in the fixed field $L_5^H$ where $H$ denotes the decomposition group
of the prime in $L_5/k$.
Now, $\p$ and $T^2+T+1$ are decomposed in $L/k$ so they are inert in 
every $R_j$, $1\leq j\leq 4$.
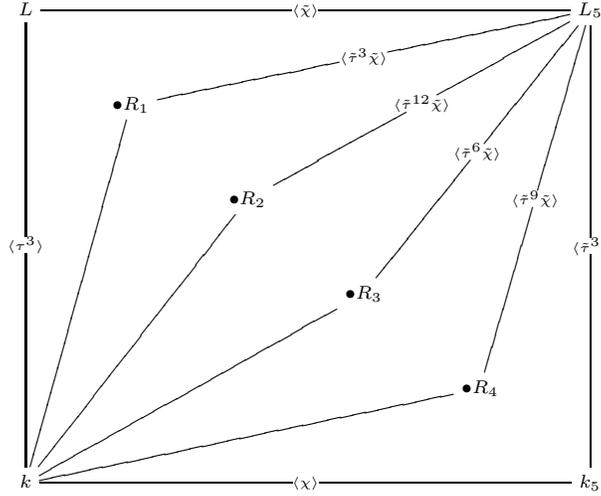
\begin{scriptsize}
\[
\xymatrix{L\ar@{-}[ddddd]|-{\langle\tau^3\rangle}
\ar@{-}[rrrrr]|-{\langle\tilde{\chi}\rangle}
&&&&&L_5\ar@{-}[ddddd]|-{\langle\tilde{\tau}^3\rangle}
\ar@{-}[lllld]|-{\langle \tilde{\tau}^3\tilde{\chi}\rangle}
\ar@{-}[llldd]|-{\langle \tilde{\tau}^{12}\tilde{\chi}\rangle}
\ar@{-}[llddd]|-{\langle \tilde{\tau}^6\tilde{\chi}\rangle}
\ar@{-}[ldddd]|-{\langle \tilde{\tau}^9\tilde{\chi}\rangle}\\
&\bullet R_1\ar@{-}[ddddl]\\
&&\bullet R_2\ar@{-}[dddll]\\
&&&\bullet R_3\ar@{-}[ddlll]\\
&&&&\bullet R_4\ar@{-}[dllll]\\
k\ar@{-}[rrrrr]|-{\langle\chi\rangle}&&&&&k_5}
\]
\end{scriptsize}

Next we compute the decomposition group 
${\mathcal D}_P$ in $k(\Lambda_M)
k_5/k$ for $P\in\{T,T+1,T^3+T^2+1,T^3+T+1, T^4+T^3+1,
T^4+T^3+T^2+T+1\}$ using formula (\ref{Eq1}). We denote by $\xi_P$
the Frobenius of $P$ in $k_5/k$ and by $\varphi_P$ the Frobenius
of $P$ in $k(\Lambda_M)/k$ (see (\ref{Eq0})).
Therefore, the Frobenius $\theta_P$
of $P$ in $k(\Lambda_M) k_5/k$ is given by $\theta_P=
\tilde{\varphi}_P \tilde{\xi}_P$. 

From (\ref{Eq2}) and from the fact that the Frobenius $\tilde{\xi}_P$
of any $P$ of degree $i$ in $k_5/k$ corresponds to $\langle \chi^i
\rangle$, we obtain $\theta_P$ and the decomposition group
${\mathcal D}_P=\langle \theta_P\rangle$ for each $P$ in
$k(\Lambda_M)k_5$ as follows
\begin{gather*}
\theta_T=\tilde{\tau}\tilde{\chi},\quad \theta_{T+1}=\tilde{\tau}^4
\tilde{\chi},\quad \theta_{T^3+T^2+1}=
\tilde{\tau}^{13}\tilde{\chi}^{3},\\
\theta_{T^3+T+1}=\tilde{\tau}^{7}\tilde{\chi}^{3},\quad
\theta_{T^4+T^3+1}=\tilde{\tau}^{9}\tilde{\chi}^{4},\quad
\theta_{T^4+T^3+T^2+T+1}=\tilde{\tau}^{6}\tilde{\chi}^{4}.
\end{gather*}

Now let $H_P$ be the subgroup of ${\mathcal D}_P$ of order $5$.
We obtain
\begin{gather}
H_T=H_{T^3+T^2+1}=H_{T^4+T^3+1}=\langle
\tilde{\tau}^{6}\tilde{\chi}\rangle,\label{Eq3}\\
H_{T+1}=H_{T^3+T+1}=H_{T^4+T^3+T^2+T+1}=\langle
\tilde{\tau}^{9}\tilde{\chi}\rangle,\nonumber
\end{gather}
and note that $\langle\tilde{\tau}^{6}\tilde{\chi}\rangle
\neq \langle \tilde{\tau}^{9}\tilde{\chi}\rangle$.

Let $R_3=L_5^{\langle \tilde{\tau}^{6}\tilde{\chi}\rangle}$ and
$R_4=L_5^{\langle \tilde{\tau}^{9}\tilde{\chi}\rangle}$.
From (\ref{Eq3}) we have that in $R_3/k$,
$T, T^3+T^2+1$ and $T^4+T^3+1$ split and in $R_4/k$,
$T+1,T^3+T+1$ and $T^4+T^3+T^2+T+1$ split. Therefore all
the primes of degree $i$ with $1\leq i\leq 4$ other than ${\eu m}$
are inert in $R_1/k$ and in $R_2/k$, where $R_1=L_5^{\langle
\tilde{\tau}^{3}\tilde{\chi}\rangle}$ and $R_2=L_5^{\langle
\tilde{\tau}^{12}\tilde{\chi}\rangle}$ and we have $\langle
\tilde{\tau}^{3}\tilde{\chi}\rangle\neq \langle
\tilde{\tau}^{12}\tilde{\chi}\rangle$. The fields $R_1$ and $R_2$
are of genus four and class number one. 

Finally, we will prove that $R_1\cong R_2$. Let $\sigma\colon k
\to k$ be given by $\sigma(T)=\frac{1}{T}$ and extend $\sigma$
to $\tilde{\sigma}\colon R_1\to \bar{k}$. Since $R_1$ and $R_2$
are the only subfields of $L_5$
of genus four and class number one, necessarily
we have $\sigma(R_1)=R_1$ or $R_2$. Now consider the primes
$T^7+T^4+1$ and $T^7+T^3+1$. Since $T^7+T^4+1\equiv
T^3+1\bmod M$ we have that $H_{T^7+T^4+1}=
\langle \tilde{\tau}^{12}\tilde{\chi}\rangle$. Therefore $T^7+T^4+1$
splits in $R_2$ and is inert in $R_1$.

Now consider the prime $T^7+T^3+1$. Since $T^7+T^3+1
\equiv T\bmod M$, it follows that $H_{T^7+T^3+1}=
\langle \tilde{\tau}^{3}\tilde{\chi}\rangle$ so that $T^7+T^3+1$
splits in $R_1$ and is inert in $R_2$. Since $\sigma(T^7+T^4+1)
=\frac{T^7+T^3+1}{T^7}$ it follows that $\sigma(R_1)=R_2$ and
$R_1\cong R_2$. This proves Theorem \ref{T2.3}. \fin

\begin{remark}\label{R2.4} {\rm{The fields $R_1$ and $R_2$
are the fields described by C. Stirpe in \cite{Sti2013}.}}
\end{remark}

The main result of this paper is a consequence of 
Theorem \ref{T2.1}, Remark \ref{R2.2} and Theorem
\ref{T2.3}.

\begin{theorem}\label{T2.4}
Up to isomorphism, there exists
exactly one function field over the
finite field of two elements of class number one and genus four. \fin
\end{theorem}

\begin{remark}\label{R2.5}{\rm{
The field $K$ (equal to either $R_1$ or $R_2$) satisfies that
${\mathcal G}=\Aut_{{\ma F}_2} K=\Aut_k K=
\Gal(K/k) \cong C_5$. Indeed, if
$|{\mathcal G}|>5$, there would exist an element of order $2$ or $3$
in ${\mathcal G}$ and if $S$ were the group generated by this
element, we would have $1<[K:K^S]=|S|<5$, thus $K^S$ 
would be of genus $0$ (see Remark \ref{R0.-1}). 
This contradicts that five is the minimum
degree of a proper subfield of $K$.
}}
\end{remark}

\end{document}